\documentclass[12pt,reqno]{amsart}

\usepackage{amsfonts}
\usepackage{amssymb}
\usepackage{amsthm}
\usepackage{amsmath}

\newtheorem{theo}{Theorem}[section]
\newtheorem{propo}[theo]{Proposition}

\newtheorem{coro}[theo]{Corollary}

\theoremstyle{definition}
\newtheorem*{defin}{Definition}

\theoremstyle{remark}
\newtheorem*{rema}{Remark}

\numberwithin{equation}{section}

\DeclareMathOperator{\GF}{GF}

\DeclareMathOperator{\diag}{diag}

\newcommand{\field}[1]{\mathbb{#1}}

\begin{document}

\title[Switching operations for Hadamard matrices]
{Switching operations for Hadamard matrices}
\author[W. P. Orrick]{W. P. Orrick}
\address{Department of Mathematics, Indiana University,
Bloomington IN 47405, USA}

\date{\today}

\begin{abstract}
We define several operations that switch substructures of Hadamard matrices
thereby producing new, generally inequivalent, Hadamard matrices.  These
operations have application to the enumeration and classification of
Hadamard matrices.  To illustrate their power, we use them to greatly
improve the lower bounds on the number of equivalence classes of Hadamard
matrices in orders 32 and 36 to 3,578,006 and 18,292,717.
\end{abstract}

\maketitle

\section{Introduction}\label{sect:intro}

Two matrices, $A$ and $B$, with entries in
the set $\{-1,1\}$ are {\em Hadamard equivalent} if $B$ can be obtained from
$A$ by some sequence of
\begin{itemize}
\item row negations,
\item column negations,
\item row permutations, and
\item column permutations.
\end{itemize}
Hadamard equivalence is so named because of its connection with Hadamard
matrices, defined as square matrices with elements equal to $\pm1$
whose rows are mutually orthogonal.  The listed moves all preserve
the property of being a Hadamard matrix.

In this paper, we describe some additional moves, called {\em switching operations,}
that preserve the
property of being a Hadamard matrix.  These operations, when applied over and over again
to a seed matrix, generally produce many inequivalent Hadamard matrices.

Furthermore, adjoining the new operations to the list above gives  new notions
of equivalence.  These weaker notions of equivalence may be useful
in the classification of Hadamard matrices  since they
partition the set of Hadamard matrices into a much smaller number of equivalence classes than does Hadamard
equivalence, but at the same time provide an effective method for
enumerating the elements of these newly defined equivalence classes.

Extensive calculation indicates that the number of Hadamard
equivalence classes that can be constructed using the new operations is enormous.
This is a big step forward since, although complete enumerations up to order 28
suggest that the number of equivalence classes grows rapidly in higher order, 
up till now there has been no general method for producing the vast numbers of
equivalence classes that we expect to exist.  The
many known Hadamard matrix construction techniques typically apply only in
scattered orders, or tend to produce Hadamard matrices with special
features such as large automorphism groups, large Hadamard submatrices,
or self-duality.

The most prolific method for constructing Hadamard
matrices has been to use two Hadamard matrices of size $n$, $A$ and $B$, to build
Hadamard matrices of size $2n$
\begin{equation}\label{kronecker}
H=\begin{bmatrix}A & PB\\A & -PB\end{bmatrix} \quad\text{and}\quad
{\widetilde H}=\begin{bmatrix}A & A\\ BP & -BP\end{bmatrix}
\end{equation}
where $P$ is any permutation matrix.  Both $A$ and $B$ can be taken from any equivalence class.
In order 32, Lin, Wallis, and Lie~\cite{LinWalLie93} produced at least 66099 inequivalent matrices
from the five equivalence classes in order 16.
Since the resulting matrices contain Hadamard submatrices of order $16$, however, they cannot
be considered generic.
In contrast, the new operations produce at least $3.57$ million equivalence classes,
most of which do not contain Hadamard submatrices of order 16.

Lam, Lam, and Tonchev have exercised great ingenuity in deriving lower
bounds on the number of Hadamard matrices of size $2n$ of the form~(\ref{kronecker}),
and have produced spectacularly large bounds in orders $40$ and higher~\cite{LamLamTon00,LamLamTon01}.  If the lessons learned from order 32 are any guide,
the true numbers of Hadamard equivalence classes in these orders are far greater still.

Our results are even more striking in orders congruent to $4\pmod{8}$ since the construction~(\ref{kronecker}) does not apply.  The previously known equivalence classes in order
36 numbered in the hundreds.  By the new methods, at least $18.29$ million
classes can be produced.

The seed matrices used to obtain all these new equivalence classes were derived from the Hadamard matrix literature up to 2005.  After this work was substantially complete, Bouyukliev, Fack, and Winne announced the classifications
of $2$-$(31,15,7)$ and $2$-$(35,17,8)$ designs with automorphisms of odd prime order.  From these
designs, they found tens of thousands of new Hadamard equivalence classes in orders 32 and 36~\cite{BouFacWin05,BouFacWin06}.  Most of these matrices have not yet been analyzed by our method.  Compared with the analysis of the dozens of previously known Hadamard equivalence classes in order 32 (excluding the matrices from construction~(\ref{kronecker})), and the hundreds of previously known H-classes in order 36, analyzing these new matrices is a major undertaking, and will require considerable optimization of our methods.  Therefore, with one important
exception, we have not used the matrices of Bouyukliev, Fack, and Winne in our enumeration, although we make a few remarks on our preliminary analysis in the next paragraph.   The exception is a matrix of order 36 in Smith class 16 (defined in Section~\ref{subsection:smith}), no previous example of which appears to have been known.  This was used as a seed matrix to produce a new family containing at least five million Hadamard equivalence classes.

In orders 4, 8, 12, 16, 20, 24,
28 the numbers of Hadamard equivalence classes are known to be 1, 1, 1, 5, 3, 60, 487~\cite{Hal61,Hal65,ItoLeoLon81,Kim89,Kim94b}.  We define a
weaker notion of equivalence, which we call {\em Q-equivalence}, by adjoining the new operations to the
operations that define Hadamard equivalence.  The numbers of Q-equivalence
classes are  1, 1, 1, 1, 1, 2, 2.  In order 32, we find that the 3.57 million known Hadamard equivalence classes are grouped into 11 Q-equivalence classes, and that in order 36, the $18.29$ million
known equivalence classes are grouped into 21 Q-equivalence classes.  As mentioned above, these numbers do not include Hadamard equivalence classes or Q-equivalence classes derived from the recently discovered matrices of Bouyukliev, Fack, and Winne.
An analysis of their matrices should provide a good test of the ideas of this  paper regarding using Q-equivalence
in classifying Hadamard matrices.  Preliminary analysis of a sampling of the the new matrices does
not turn up any new large Q-classes, but does indicate the presence of a large number of
new small Q-classes (perhaps in the hundreds or more).  We intend to make a complete enumeration of these, and a full analysis of all the new matrices.  The results will be presented in a follow-up to the present paper.

\section{Overview of switching}
Suppose that an $n\times n$ Hadamard matrix can be put in the form
\setcounter{MaxMatrixCols}{15}
\begin{equation}\label{ccForm}
\begin{bmatrix}
\mathbf{1} & \cdots & \mathbf{1} & - & \cdots & - & - & \cdots & - & 1 & \cdots & 1 \\
\mathbf{1} & \cdots & \mathbf{1} & - & \cdots & - & 1 & \cdots & 1 & - & \cdots & - \\
\mathbf{1} & \cdots & \mathbf{1} & 1 & \cdots & 1 & - & \cdots & - & - & \cdots & - \\
\mathbf{1} & \cdots & \mathbf{1} & 1 & \cdots & 1 & 1 & \cdots & 1 & 1 & \cdots & 1 \\
 & a_5  & & & b_5 & & & c_5 & & & d_5\\
 & \vdots & & & \vdots & & & \vdots & & & \vdots \\
 & a_n & & & b_n & & & c_n & & & d_n\\
\end{bmatrix}
\end{equation}
where $a_i$, $b_i$, $c_i$, and $d_i$ are $\{-1,1\}$-vectors of length $n/4$.  The columns of
the matrix
have been grouped into four sets of $n/4$ columns each.  A new,
generally inequivalent, Hadamard matrix can be obtained by
negating the $4\times \frac{n}{4}$ block of 1s in the upper left corner (shown in boldface).
We call this operation {\em switching a closed quadruple.}
\setcounter{MaxMatrixCols}{10}

Suppose instead that we can put the matrix in the form of Figure~\ref{figure:hallSet}
where the $A_{ij}$ are square matrices of size $(n-4)/4$.  A new, often inequivalent, matrix can be obtained by
negating the all 1 block of size $4\times\frac{n-4}{4}$ contained in the first four rows, and
the all 1 block of size $\frac{n-4}{4}\times4$ contained in the first four columns (both shown
in boldface).  We call this operation {\em switching a Hall set.}
\begin{figure}
\setcounter{MaxMatrixCols}{16}
\begin{equation}
\begin{bmatrix}
1 & - & - & - & \mathbf{1} & \cdots & \mathbf{1} & 1 & \cdots & 1 & 1 & \cdots & 1 & - & \cdots & - \\
- & 1 & - & - & \mathbf{1} & \cdots & \mathbf{1} & - & \cdots & - & - & \cdots & - & - & \cdots & - \\
- & - & 1 & - &\mathbf{1} & \cdots & \mathbf{1} & - & \cdots & - & 1 & \cdots & 1 & 1 & \cdots & 1 \\
- & - & - & 1 &\mathbf{1} & \cdots & \mathbf{1} & 1 & \cdots & 1 & - & \cdots & - & 1 & \cdots & 1 \\
\mathbf{1} & \mathbf{1} & \mathbf{1} & \mathbf{1} \\
\vdots & \vdots & \vdots & \vdots & & A_{11} & & & A_{12} & & & A_{13} & & & A_{14} \\
\mathbf{1} & \mathbf{1} & \mathbf{1} & \mathbf{1} \\
- & 1 & 1 & - \\
\vdots & \vdots & \vdots & \vdots & & A_{21} & & & A_{22} & & & A_{23} & & & A_{24} \\
- & 1 & 1 & - \\
- & 1 & - & 1 \\
\vdots & \vdots & \vdots & \vdots & & A_{31} & & & A_{32} & & & A_{33} & & & A_{34} \\
- & 1 & - & 1 \\
1 & 1 & - & - \\
\vdots & \vdots & \vdots & \vdots & & A_{41} & & & A_{42} & & & A_{43} & & & A_{44} \\
1 & 1 & - & - \\
\end{bmatrix}
\end{equation}
\setcounter{MaxMatrixCols}{10}

\caption{Switching a Hall set}
\label{figure:hallSet}
\end{figure}

Justification for these claims and further elaboration are given in the subsequent sections.

\section{Closed quadruples and Hall sets}
\subsection{3-normalization}
Let $H$ be a Hadamard matrix of size $n$.  Denote its rows by $h_i$ and its
elements by $h_{ij}$.  Define the Hadamard product of two vectors to be
\begin{equation*}
(a_1,\ldots,a_n)\circ (b_1,\ldots, b_n):=(a_1 b_1,\ldots,a_n b_n).
\end{equation*}
Let $j_k$ be the all 1 vector of length $k$.

\begin{defin}
A Hadamard matrix of size $n$ is
{\em 3-normalized on rows $(i,j,k)$} if, in every column $\ell$, the set
$\{h_{i\ell},h_{j\ell},h_{k\ell}\}$ contains an even number of $-1$s, or
equivalently if $h_i\circ h_j\circ h_k=j_n$.
\end{defin}

3-normalization is a normalization of the columns.  A 3-normalized matrix
remains 3-normalized if any two of the rows $i$, $j$, $k$ or of any
single row other than $i$, $j$, $k$ is negated.
Note that 3-normalization was introduced in~\cite{OrrSol07}.
The definition given here is slightly weaker in that it makes no stipulation
that the row sums be positive, and does not impose any particular ordering
on the columns.  In the next paragraph we restate some needed results from~\cite{OrrSol07}.

The {\em field structure} $(C_1,C_2,C_3,C_4)$ of a 3-normalized
Hadamard matrix of size $n$  is the partition of the
set of columns $c$ into four classes, $C_i$, accordingly as
$(h_{jc},h_{kc},h_{\ell c})=(1,1,1)$, $(-1,-1,1)$, $(-1,1,-1)$,
or $(1,-1,-1)$.  The four classes are
called {\em fields} and are all of length $n/4$.  In a row $r\notin\{j,k,\ell\}$
the sum of the elements in a field is the same
for each of the four fields in the row.  This follows from orthogonality of
row $r$ with rows $j$, $k$, $\ell$.  Since
the sum of the entries in a field is even if $n/4$ is even,
and odd if $n/4$ is odd,
the row sum of row $r\notin\{j,k,\ell\}$ must be congruent to $n\pmod{8}$.

A quadruple of rows, $(i,j,k,\ell)$ of a Hadamard matrix $H$ of size $n$ is said to
be of {\em type $r$}, $0\le r\le n/8$, if exactly $4r$ of
the entries in $h_i\circ h_j \circ h_k\circ h_\ell$ equal $-1$
or exactly $4r$ entries equal $+1$.  This notion was introduced by Kimura~\cite{Kim94b}.

\begin{defin}
A quadruple of rows, $(i,j,k,\ell)$ of a Hadamard matrix $H$ of size $n$, is {\em closed} if
$h_i\circ h_j \circ
h_k\circ h_\ell=\pm j_n$.
\end{defin}

A closed quadruple is a quadruple of type 0.
Thus if $H$ is 3-normalized on three rows of a closed quadruple, then the
fourth will consist entirely of 1s or entirely of $-1$s.  The field structure is
independent of which three rows of the closed quadruple are chosen.

Quadruples of type 1 will also play an important role in what follows.  They were used extensively by Hall in the classification of Hadamard matrices of order 20~\cite{Hal65} and by Kimura in the classification for order 28~\cite{Kim94a,Kim94b}.  If $H$ is 3-normalized on three rows of a type-1 quadruple, then the fourth row will contain one odd-sign entry in each of the fields induced
by the 3-normalization.  Kimura and Ohmori referred to such quadruples
as {\em Hall sets}~\cite{KimOhm86}.

\begin{propo}
If a Hadamard matrix of size $n$ has a closed quadruple, then $n=4$ or
$n\equiv0\pmod{8}$.
\end{propo}
\begin{proof}
Let $(i,j,k,\ell)$
be the closed quadruple.  3-normalize the matrix
on rows $i$, $j$, $k$ so that
$h_\ell=\pm j_n$.  Orthogonality implies that all rows except for $h_\ell$
have row sum 0.
All row sums of rows other than $i$, $j$, $k$ must be congruent
to $n\pmod{8}$.  If $n>4$ this can only happen when $n\equiv0\pmod{8}$.
\end{proof}

\subsection{Obtaining new Hadamard matrices by switching closed quadruples}

\begin{defin}
Let $H$ be a Hadamard matrix of size $n$ which has a closed quadruple, $Q$.
Let $(C_1,C_2,C_3,C_4)$ be the partition of columns induced by 3-normalization
on $Q$.
{\em Switching the closed quadruple} $Q$ means negating all the elements
$h_{rc}$, where $r\in Q$ and $c\in C_i$ for some
$i\in\{1,2,3,4\}$.
\end{defin}

\begin{propo}\label{propo:closed}
The matrix produced by switching a closed quadruple $Q$ in a
Hadamard matrix $H$ is a Hadamard matrix.
\end{propo}
\begin{proof}
Any matrix containing a closed quadruple is Hadamard equivalent to one of the form~(\ref{ccForm}).  It is evident that switching preserves orthogonality of the columns in that matrix.  Since column orthogonality is preserved under the operations needed to put $H$ in the form~(\ref{ccForm}), the conclusion holds generally.
\end{proof}

It appears that when $n>8$, switching always produces a Hadamard
matrix that is inequivalent to the original Hadamard matrix.

Note that the equivalence class of the Hadamard matrix produced by
switching $Q$ is independent of which of the four fields
$C_i$ we choose to negate.  To see this, note that
negating the closed quadruple elements in $C_2$ is equivalent
to first negating the closed quadruple elements in $C_1$, then
negating all four rows of the
closed quadruple, and finally performing
a certain permutation of the rows of $Q$.  The
same holds for $C_3$ and $C_4$.

\subsection{More general row switching operations}
\label{subsection:general}
It was observed by Denniston~\cite{Den82} in connection with symmetric $(25,9,3)$ designs that,
starting from a design, a new inequivalent design
can be obtained by switching a substructure known as an oval.  
Denniston's switching operation can be thought of as an operation that permutes certain elements
of the incidence matrix of the design.  In fact, if zeroes are replaced with $-1$s in the incidence
matrix of a $(25,9,3)$ design, ovals satisfy our definition of a closed quadruple, and our switching
operation is equivalent to Denniston's.

We can formulate more general switching operations acting on more general structures, which
we will refer to generically as ``designs.''  Consider a set, $\mathcal{M}$, of matrices of a fixed size which
represent the designs in question.  If $R\in\mathcal{M}$, we suppose that the elements of $R$ 
are taken from some set $\mathcal{S}$, that the rows satisfy some set of properties
$\mathcal{P}$, and furthermore, that $R$ satisfies
\begin{equation*}
R^{\rm{T}}R=M
\end{equation*}
where $M$ is some fixed matrix.  For example, if $\mathcal{M}$ represents $(25,9,3)$ designs,
then $\mathcal{S}=\{0,1\}$, the set $\mathcal{P}$ contains the property that every row of $R\in\mathcal{M}$
has exactly nine 1s, and $M=6I+3J$ where $J$ is the $25\times25$ all 1 matrix.  If $\mathcal{M}$ represents Hadamard matrices of size $n$,
then $\mathcal{S}=\{-1,1\}$, the set $\mathcal{P}$ is empty, and $M=nI$.  Many other types of
matrices and designs, including certain D-optimal designs
can also be defined within this framework.

Let $R\in\mathcal{M}$ and partition the incidence matrix into two submatrices, $A$ and $X$,
\begin{equation*}
R=\begin{bmatrix} A\\ X\end{bmatrix}.
\end{equation*}
Now suppose that $B$ is a matrix of the same dimensions as $A$, with elements taken from the
same set $\mathcal{S}$, satisfying the same
properties $\mathcal P$, and that $B^{\rm{T}}B=A^{\rm{T}}A$.  Then the matrix obtained
from $R$ by replacing $A$ with $B$ is also a matrix of the original type.

Suppose for example, that $R$ is an $n\times n$ Hadamard matrix with an $m\times n$ submatrix $A$ whose columns are all identical to columns of a particular $m\times m$ Hadamard matrix $H_m$
or to negations of such columns.  Let $H'_m$ be another $m\times m$ Hadamard matrix.  Denote
column $j$ of $H_m$ by $v_j$ and column $j$ of $H'_m$ by $v'_j$.  Define $[a,b]$ to be the set of
integers $i$ satisfying $a\le i\le b$.  For $j=1,\ldots,n$, let column $j$ of $A$ be $\sigma(j) v_{a(j)}$ where $\sigma:[1,n]\rightarrow\{-1,1\}$ and $a:[1,n]\rightarrow[1,m]$.  Let $B$ be the matrix with
the same dimensions as $A$ whose columns are $\sigma(j)v'_{a(j)}$ for $j=1,\ldots,n$.  Then $B$ will
satisfy $B^{\rm{T}}B=A^{\rm{T}}A$, and so we may use it to obtain a new Hadamard matrix
of order $n$.

Note that if $m=1$ and we let $H_1=\begin{bmatrix}1\end{bmatrix}$ and $H'_1=\begin{bmatrix}-1\end{bmatrix}$ then the above operation amounts to negation
of a row.  Likewise, if $m=2$ and $H_2=\begin{bmatrix}1 & 1\\1 & -1\end{bmatrix}$ while
$H'_2=\begin{bmatrix}1 & -1\\1 & 1\end{bmatrix}$ then the above operation amounts to
swapping two rows.

Switching closed quadruples is an instance of the $m=4$ case.
Let $A$ be a $4\times n$
matrix whose columns, or their negations, are columns of $H_4$, a $4\times4$ Hadamard matrix.
Orthogonality of the rows of $A$ implies that $A$ is a closed quadruple.
Now negating one column of $H_4$ and using the resulting matrix to
construct $B$, has the effect of switching the closed quadruple formed by
the rows of $A$.  Thus, in some sense switching closed quadruples is a natural extension
of the operations of row negation and row permutation.

\subsection{Closed quadruples and Hadamard submatrices}
\label{subsection:tensor}
There is an additional sense in which switching closed quadruples is a natural
extension of the operation of row permutation.
Consider the matrix $H$ of size $2n$ defined in equation~(\ref{kronecker}).
One may negate or permute the columns of $A$ or $B$
without changing the equivalence class of $H$.  One may also negate a row
of $PB$ (or of $A$) without changing the equivalence class of $H$.  The reason is that
negating row $j$ of $PB$
amounts to swapping rows $j$ and $n+j$ of $H$.

On the other hand, changing the permutation $P$, for example by performing the additional
row swap $(i,j)$, usually
does change the equivalence class of $H$.  The additional swap will affect
four rows of $H$, namely $i$, $j$, $i+n$, $j+n$.  These four rows
form a closed quadruple.  One of the four fields of this quadruple
is the set of columns of $H$ in which rows $i$ and $j$ of $PB$
differ.  We make the switch that negates the entries in rows
$i$, $j$, $i+n$, and $j+n$ that lie within this field.  The result is identical
to the result of swapping rows $i$ and $j$ of $PB$.
Therefore, in this context, switching a closed quadruple
amounts to swapping a pair of rows in one of the two matrices
from which $H$ was constructed.

\subsection{Properties of Hall sets}
Hall sets play the role for matrices of order $n\equiv4\pmod{8}$
that closed quadruples play for matrices of order $n\equiv0\pmod{8}$.

Hall sets can be found both in Hadamard matrices of order
$n\equiv0\pmod{8}$ and in those of order $n\equiv4\pmod{8}$.
Four columns are singled out in the definition of a Hall set, namely the
columns whose sign in the Hadamard product differs from the sign of all the other columns.  When $n\equiv4\pmod{8}$
these form a Hall set in the columns, as shown
by Kimura and Ohmori~\cite{KimOhm86}.  For convenience of the reader, we
reprove this here.  We include the corresponding result for $n\equiv0\pmod{8}$
for good measure.  Define the {\em Hall columns} to be the four distinguished columns.  There is one Hall
column in each field.

\begin{propo}
Let $H$ be a Hadamard matrix of order $n$.  If $n\equiv0\pmod{8}$ then the
Hall columns form a closed quadruple.  If $n\equiv4\pmod{8}$ then the
Hall columns form a Hall set.
\end{propo}

\begin{proof}
We assume without loss of generality that $H$ is 3-normalized on three rows of the Hall set.
Consider a row not contained in the Hall set.   Let $x_i$ denote the element of that row in the Hall
column of field $i$.  Let $a_i$ denote the sum of the remaining elements of field $i$.
Then orthogonality with the Hall set rows implies
\begin{align}
&x_1+a_1=x_2+a_2=x_3+a_3=x_4+a_4\\
& a_1+a_2+a_3+a_4=x_1+x_2+x_3+x_4,
\end{align}
which implies that the row sum, which must be congruent to $n\pmod{8}$, equals $2(x_1+x_2+x_3+x_4)$.  Hence the product $x_1x_2x_3x_4$
is positive for $n\equiv0\pmod{8}$ and negative for $n\equiv4\pmod{8}$.  In each row of the Hall set, the
product of the four elements in Hall columns is always positive, so the result follows.
\end{proof}

\begin{rema}
When $n\equiv0\pmod{8}$ the existence of a Hall set implies the existence of a closed quadruple
in the columns, but the converse is not true.  The existence of a closed quadruple does not imply
the existence of a corresponding Hall set. \footnote{Sylvester matrices exhibit this in extreme form.
One can show by induction on $k$ that the Sylvester matrix of size $2^k$  has $\frac{1}{4}\binom{2^k}{3}$
closed quadruples, and that any of its other row quadruples is of type $2^{k-3}$ (which means
its Hadamard product  has as many
entries $+1$ as $-1$).  Therefore, if $k\ge4$, the Sylvester matrix has many closed quadruples, but
no Hall sets.  Since Sylvester matrices are self dual, the same is true of column quadruples.}
\end{rema}

Henceforth we will consider the $n\equiv4\pmod{8}$ case, and when we speak of a Hall set, we
will mean both the four rows of the set and the four corresponding Hall columns.

By permuting the Hall rows and columns to the top- and leftmost positions
and normalizing appropriately we obtain the form
\begin{equation}\label{hallSet}
H=\begin{bmatrix}
H_4 & F_1 & F_2 & F_3 & F_4\\
G_1 & A_{11} & A_{12} & A_{13} & A_{14} \\
G_2 & A_{21} & A_{22} & A_{23} & A_{24} \\
G_3 & A_{31} & A_{32} & A_{33} & A_{34} \\
G_4 & A_{41} & A_{42} & A_{43} & A_{44}
\end{bmatrix}
\end{equation}
where
\begin{align}\label{Fdefs}
H_4=\begin{bmatrix}1 & - & - & -\\ - & 1 & - & -\\ - & - & 1 & -\\ - & - & - & 1\end{bmatrix} \qquad
F_1&=\begin{bmatrix}1 & \ldots & 1\\1 & \ldots & 1\\1 & \ldots & 1\\1 & \ldots & 1\end{bmatrix}\qquad
&F_2&=\begin{bmatrix}1 & \ldots & 1\\- & \ldots & -\\- & \ldots & -\\1 & \ldots & 1\end{bmatrix}\notag\\
F_3&=\begin{bmatrix}1 & \ldots & 1\\- & \ldots & -\\1 & \ldots & 1\\- & \ldots & -\end{bmatrix}\qquad
&F_4&=\begin{bmatrix}- & \ldots & - \\- & \ldots & - \\1 & \ldots & 1 \\1 & \ldots & 1\end{bmatrix},
\end{align}
$G_1=F_1^{\mathrm{T}}$, $G_j=-F_j^{\mathrm{T}}$ for $j\in\{2,3,4\}$, 
and $A_{ij}$ are submatrices whose row and column sums equal $2$ when
$i=j$ and $0$ when $i\ne j$.

\begin{defin}
By {\em switching a Hall set} in the matrix $H$ defined in eqn.~(\ref{hallSet}) we mean
the operation of replacing $F_i$ by its negation and $G_i$ by its negation for
one of the choices $i=1,2,3,4$.
\end{defin}

The four possible negations of the definition produce equivalent matrices.  The proof of this is similar
to the proof of the analogous property of closed quadruples given in the discussion
following Proposition~(\ref{propo:closed}).  Switching is well defined even when the
Hall rows and columns do not appear in positions 1--4 or when the normalization is different
from the one in~(\ref{hallSet}).  We need only apply a signed permutation to put the matrix
into the form~(\ref{hallSet}), switch as in the definition, and then apply
the inverse signed permutation.

\begin{propo}\label{propo:hallSwitch}
The matrix produced by switching a Hall set in a
Hadamard matrix is a Hadamard matrix.
\end{propo}

\begin{proof}
We will assume the form~(\ref{hallSet}) since the conclusion is unaffected by the permutations and
negations needed to convert the matrix to that form.  When $i\ne j$,  the rows of $F_j$
are orthogonal to the rows of $A_{ij}$ as the latter have row sum $0$.  Therefore, negating $F_j$ does not alter the orthogonality of rows 1--4 of $H$ with the rows of $H$ contained in the block $\begin{bmatrix}G_i & A_{i1} & A_{i2} & A_{i3} & A_{i4}\end{bmatrix}$. Row $k$ ($k=1,2,3,4$) of $F_j$ has inner product $\pm2$ with any of the rows of $A_{jj}$ while row $k$ of
$H_4$ has inner product $\mp2$ with any of the rows of $G_j$.  Negating both $F_j$ and $G_j$
produces sign changes in these inner products that produce opposite contributions to any of the inner products
of rows 1--4 of $H$ with the rows of $H$ contained in the block $\begin{bmatrix}G_j & A_{j1} & A_{j2} & A_{j3} & A_{j4}\end{bmatrix}$.
\end{proof}

Examples are known where switching a Hall set in a Hadamard matrix $H$ produces
a matrix equivalent to $H$.  In general, however, one obtains an inequivalent matrix.

\section{Invariants}
First we prove that the number of closed quadruples in a Hadamard matrix of size $16k+8$ is invariant
under switching closed quadruples.  Second,
we show that the the binary, doubly even, self-dual code associated to the transpose of a Hadamard matrix
of size $16k+8$ is unchanged by switching a closed quadruple of that matrix.
Finally  we show that  the integer equivalence class is preserved under switching of Hall sets.

\subsection{A closed quadruple switching invariant for $n\equiv8\pmod{16}$}
We will need to understand the ways that
closed row quadruples may overlap within a Hadamard matrix.
\begin{propo}
Suppose $(i,j,k,\ell)$ and $(i',j',k',\ell')$ are distinct closed quadruples
with nonempty intersection.  Then the number of rows common to the two
quadruples is 2 if $n\equiv8\pmod{16}$ and 1 or 2 if $n\equiv0\pmod{16}$.
\end{propo}
\begin{proof}
The number of common rows cannot be 3 since the fourth row of a closed
quadruple is determined, up to sign, by the other three, and the two
quadruples are assumed distinct.  Therefore the number of common rows
must be either 1 or 2.

We will show that if the number of common rows is 1, then $n\equiv0\pmod{16}$.

Assume the number of common rows to be 1 and let $n=8r$.
Take the two quadruples to be $(1,2,3,4)$
and $(1,5,6,7)$, and 3-normalize the matrix on rows 2, 3, 4.  Normalize row 1
to have positive entries.  By suitable column permutations, the structure of
the first five rows can be brought to the form:

\begin{center}
\begin{tabular}{lcrrcrrcrrcrr}
1. && $1_r$ & $1_r$ && $1_r$ & $1_r$ && $1_r$ & $1_r$ && $1_r$ & $1_r$ \\
2. && $1_r$ & $1_r$ && $-1_r$ & $-1_r$ && $-1_r$ & $-1_r$ && $1_r$ & $1_r$ \\
3. && $1_r$ & $1_r$ && $-1_r$ & $-1_r$ && $1_r$ & $1_r$ && $-1_r$ & $-1_r$ \\
4. && $1_r$ & $1_r$ && $1_r$ & $1_r$ && $-1_r$ & $-1_r$ && $-1_r$ & $-1_r$ \\
5. && $1_r$ & $-1_r$ && $1_r$ & $-1_r$ && $1_r$ & $-1_r$ && $1_r$ & $-1_r$ \\
\end{tabular}
\end{center}

The form of row 5 is a consequence of the fact that the sum of elements in
each of the four fields must be zero.  Since $(1,5,6,7)$ is closed,
the Hadamard product of rows 6 and 7 equals either row 5 or its negation.
By normalizing row 7 appropriately we may assume the former.
Consider the two subfields that compose the first field
in the above structure.  They will be further subdivided as

\begin{center}
\begin{tabular}{lcrrcrrc}
5. && $1_a$ & $1_{r-a}$ && $-1_b$ & $-1_{r-b}$ &\ldots\\
6. && $1_a$ & $-1_{r-a}$ && $1_b$ & $-1_{r-b}$ &\ldots\\
7. && $1_a$ & $-1_{r-a}$ && $-1_b$ & $1_{r-b}$ &\ldots\\
\end{tabular}
\end{center}
The subfields composing the remaining three fields will be subdivided
similarly.  Because there
are $r$ 1s per field in rows 6 and 7, just as in row 5, we have the
constraints $a+b=r$ and $a+(r-b)=r$.  Therefore $a=b=r-a=r-b=r/2$ and
hence $r$ is even.  Consequently $n\equiv0\pmod{16}$.
\end{proof}

Note that all of the degrees of overlap between closed quadruples
allowed by the Proposition occur in practice.

\begin{propo}\label{propo:CQinvariant}
Let $n\equiv8\pmod{16}$.  Let $H$ be a Hadamard matrix of size $n$ which has
a closed row quadruple $Q$.  Switching $Q$ does not change the number
of closed row quadruples in $H$.
\end{propo}
\begin{proof}
In the matrix obtained from $H$ by switching $Q$, the rows of $Q$ still form
a closed quadruple.
Also, any quadruple, whether closed or not, that doesn't involve any rows of $Q$
is unaffected by switching.  The only way the number of closed quadruples
could change is if a closed quadruple were created or destroyed by switching $Q$.
Such a closed quadruple would have to overlap $Q$ (either before or after switching)
and would therefore
share exactly two of $Q$'s rows.  However, the Hadamard product of any pair
of rows in $Q$ is not altered by negation of any of the fields of $Q$.  Hence the
Hadamard product of the four rows
of a putative overlapping quadruple would be unchanged by such a negation.
Therefore, any closed quadruple overlapping $Q$ in two rows remains closed
after switching $Q$.  Likewise, any quadruple overlapping $Q$ in two rows
which is not closed initially, will not be closed after switching $Q$.
\end{proof}

It is worth pointing out that switching a closed {\em column} quadruple does change
the number of closed row quadruples in general.  Furthermore, switching closed row
quadruples generally does change the number of closed row quadruples when
$n\equiv0\pmod{16}$.  For example, when $n=16$, the five equivalence classes
of Hadamard matrices have 140, 76, 44, 28, and 28 closed row quadruples.  Each
of these five classes can be obtained starting from any of the others and performing
a series of switches of closed row quadruples.

\subsection{Invariant codes}
\label{subsection:codes}
Codes can be associated with Hadamard matrices, and are useful in their classification.  
For our purposes, codes can be thought of as collections of vectors over some finite field. The vectors in a code are
called {\em codewords,} and the {\em weight} of a codeword is the number of its entries that
are non-zero.    The {\em support} of a codeword is the set of positions in which it
has a non-zero entry.  

One way to associate a {\em linear} code with a Hadamard matrix of size $n$ is to normalize the columns of the matrix so
that all entries in the first row equal 1, then to change all $-1$ entries to 0, and finally to take the linear
span of the rows of the resulting matrix over some finite field $\field{F}_p$ where $p$ is a prime.
One could equally well normalize on a row other than the first and the resulting code would
be the same.
The {\em dimension} of such a linear code is its dimension as a subspace of $\field{F}_p^n$.  If $\mathcal{C}$ is such a linear code, then its {\em dual code,} $\mathcal{C}^\perp$ is the subspace of $\field{F}_p^n$ consisting of all vectors orthogonal to
all codewords in $\mathcal{C}$.  Basic linear algebra implies that the dimensions of a code and its dual satisfy $\dim(\mathcal{C})+\dim(\mathcal{C}^\perp)=n$.  If $\mathcal{C}\subseteq\mathcal{C}^\perp$ then
$\mathcal{C}$ is said to be {\em self-orthogonal}.  If $\mathcal{C}=\mathcal{C}^\perp$ then
$\mathcal{C}$ is said to be {\em self-dual}.

We will only consider binary codes ($p=2$) in this paper, but it should be noted that codes over
$\field{F}_p$, $p$ an odd prime, are closely connected with integer equivalence, which is discussed
in the next section.  The {\em 2-rank} of a Hadamard matrix is the same as the dimension of its associated binary
code.  Two binary codes are {\em isomorphic} if one can be converted to the other by a 
permutation of coordinate positions.  The following result is proved (in greater generality) in many places.  (For example, see~\cite{Lan83}, Section 2.3.)

\begin{theo}\label{theo:HadamardCodes} Let $H$ be a Hadamard matrix of size $n$.  Let $\mathcal{C}$ be a binary code associated to $H$ as described above.  Then,
\begin{enumerate}
\item if $n\equiv4\pmod{8}$ then $\mathcal{C}=\{1_n\}^\perp$ which implies $\dim(\mathcal{C})=n-1$;
\item if $n\equiv0\pmod{8}$ then $\mathcal{C}$ is self-orthogonal which implies $\dim(\mathcal{C})\le n/2$;
\item if $n\equiv8\pmod{16}$ then $\mathcal{C}$ is self-dual which implies $\dim(\mathcal{C})=n/2$.
\end{enumerate}
\end{theo}

Since all Hadamard matrices of a size congruent to $4\pmod{8}$ have the same binary code, the
binary code does not help with classification (although codes over other fields may).  For the
present, we focus on binary codes associated with matrices of size $n\equiv0\pmod{8}$.   It is not hard to show that such codes are {\em doubly-even,} that is, all of their code words have weight divisible by 4.
For an illustration that various 2-ranks allowed by
Theorem~\ref{theo:HadamardCodes} do occur in practice, we consider consider some results
discussed by Assmus and Key in~\cite{AssKey92a,AssKey92b}.   They note that the five non-equivalent Hadamard matrices of size 16 have binary codes of dimensions 5, 6, 7, 8, and 8.  Only the last two are self-dual, and they turn out not to be isomorphic.   Contrast this with size 24 where the 60 non-equivalent Hadamard matrices must {\em all} have self-dual codes of dimension 12.  Assmus and Key proved that these 60 classes of matrices are associated with six different doubly-even, self-dual, binary codes.

Jennifer Key pointed out~\cite{Key05} that when $H$ is a Hadamard matrix of size $24$, the number of closed quadruples coincides with the number of code words of weight 4 in the binary code associated with the {\em columns} of $H$.  (We might also call this the code associated with $H^\mathrm{T}$.)  We elaborate a bit on her observation, which reflects a general phenomenon for matrices of size $n\equiv8\pmod{16}$.

\begin{propo}\label{propo:wordsToQuads} Let $H$ be a Hadamard matrix of size $n\equiv0\pmod{8}$.  Let $\mathcal{C}$ be
the linear binary code constructed from the columns of $H$.  That is, $\mathcal{C}$ is the linear
span over $\field{F}_2$ of the columns of a matrix $A$  formed
by normalizing $H$ so that one of its columns consists entirely of 1s and then changing $-1$s to $0$s.  Let $\{i,j,k,\ell\}$ be the support of a weight 4 codeword in
$\mathcal{C}$.  Then rows $i$, $j$, $k$, and $\ell$ of $H$ form a closed quadruple.
\end{propo}

\begin{proof} Since
$\mathcal{C}$ is self-orthogonal, the weight 4 codeword with support $\{i,j,k,\ell\}$ is orthogonal to every column of $A$.  This means
that every column of $A$ has an even number of 1s in positions $i$, $j$, $k$, and $l$, which implies
the result.
\end{proof}

This result has a partial converse with self-duality of $\mathcal{C}$ being the needed additional assumption.

\begin{propo}\label{propo:quadsToWords}
Let $H$, $A$, and $\mathcal{C}$ be defined as in Proposition~\ref{propo:wordsToQuads} and suppose in
addition that $\mathcal{C}$ is self-dual.  Let $i$, $j$, $k$, and $\ell$ label the rows of a closed
quadruple in $H$.  Then $\{i,j,k,\ell\}$ is the support of a weight 4 codeword in $\mathcal{C}$.
\end{propo}

\begin{proof} Since one column of $A$ consists entirely of 1s, and since $i$, $j$, $k$, and $\ell$ label a closed quadruple, every column of $A$ has an even number of 1s among the positions $i$, $j$, $k$, and $\ell$.
Let $\mathbf{c}$ be the vector in $\field{F}_2^n$ with support $\{i,j,k,\ell\}$.  Then $\mathbf{c}$ is orthogonal to every column of $A$ and therefore to every codeword in $\mathcal{C}$.  Hence $\mathbf{c}\in\mathcal{C}^\perp$.  Since $\mathcal{C}$ is self-dual, we also have $\mathbf{c}\in\mathcal{C}$.
\end{proof}

We have established a one-to-one correspondence between the closed quadruples of a Hadamard matrix $H$ of size $n\equiv0\pmod{8}$ and weight 4 code words in the binary code associated to the columns of $H$, provided that that code is self-dual.   This correspondence therefore holds for all Hadamard matrices of size congruent to $8\pmod{16}$.

We finally investigate the effect of switching a closed quadruple of $H$ on the binary code associated to the columns of $H$.  In this connection, we note that the closed quadruple switching operation was defined and used in the coding theory context by Phelps, Rif{\`a}, and Villanueva~\cite{PheRifVil05}.  
They were concerned with Hadamard matrices of size $2^t$ whose codes can range in dimension
from $t+1$ to $2^{t-1}$.  Starting with a code of minimal dimension, corresponding to the Sylvester
matrix, they produced codes, and the corresponding matrices, of the next two higher dimensions
by switching.  For further details, see Lemmas 4.2 and 4.3 of~\cite{PheRifVil05}.  Note that closed quadruples correspond to subcodes of dimension three.

Our focus in this paper will be on codes at the opposite end of the range of possible dimensions, that is, on the self-dual codes.  We have the following inclusion of codes:

\begin{propo}\label{propo:codeInclusion} Let $H$, $A$, and $\mathcal{C}$ be defined as in Proposition~\ref{propo:wordsToQuads} and suppose in
addition that $\mathcal{C}$ is self-dual.  Let $H'$ be a Hadamard matrix obtained from $H$ by
switching a closed quadruple, and let $\mathcal{C}'$ be the code associated to the columns of $H'$.
Then $\mathcal{C}'\subseteq\mathcal{C}$.  Furthermore, if $\mathcal{C}'\ne\mathcal{C}$, then
$\mathcal{C}$ is spanned by $\mathcal{C}'$ and a particular weight 4 vector.
\end{propo}

\begin{proof}
We may assume that $H$ has been normalized so that all entries in its first column equal 1.   The
matrix $A$ is then obtained from $H$ simply by replacing $-1$s with 0s.  The code $\mathcal{C}$ is
the span of the columns of $A$.

Let $(i,j,k,\ell)$ be a closed quadruple of $H$.  Proposition~\ref{propo:quadsToWords} asserts that $\mathcal{C}$ contains a codeword $\mathbf{c}$ with support $\{i,j,k,\ell\}$.   Let $\{C_1,C_2,C_3,C_4\}$ be the partition into fields of the set of columns of $H$ induced by the closed quadruple $(i,j,k,\ell)$.  Switching means negating
all matrix elements in rows $i$, $j$, $k$, $\ell$ and in the columns of one of the $C_m$, $m=1$, $2$, $3$, or $4$.  Should column 1 be one of the affected columns, the normalization of the resulting matrix, $H'$, will no longer be such
that column 1 contains 1s only.  To restore the normalization, we simply negate rows $i$, $j$, $k$, $\ell$.
The net result will be that all columns {\em but} those of $C_m$ are affected by the switching.

At any rate, the matrix $A'$, obtained from $H'$ by changing $-1$s to 0s, will differ from $A$ 
only in that the elements in rows $i$, $j$, $k$, $\ell$ and in a certain subset of the columns
will have been changed to their complements ($0\rightarrow1$, $1\rightarrow0$).  This change
can be effected by adding the vector $\mathbf{c}$ to the appropriate columns of $A$.  Therefore, $\mathcal{C}'$, which is the span of the columns of $A'$, is spanned by a set of linear combinations of codewords
in $\mathcal{C}$.  Hence $\mathcal{C}'\subseteq\mathcal{C}$.  Finally, $\mathcal{C}$ is
clearly the span of $\mathcal{C}'\cup\{\mathbf{c}\}$.
\end{proof}

We note that the code $\mathcal{C'}$ obtained in the above proof depends on which of the four fields, $C_m$,
was used in the switching.  Nevertheless, the isomorphism class of the code will be independent
of this choice.

As an illustration of the use of Proposition~\ref{propo:codeInclusion} consider the codes associated
with the five equivalence classes of $16\times16$ Hadamard matrices.  Matrices in either of the two classes
associated with self-dual codes have 28 closed quadruples.  Switching any of these 28 quadruples
produces a matrix in the class corresponding to the code of dimension 7.  We therefore conclude
that the code of dimension 7 is a subspace of both of the codes of dimension 8, and that each of the latter
can be obtained by augmenting the code of dimension 7 with a  weight 4 vector whose support
corresponds to a suitable closed quadruple.

A corollary of Proposition~\ref{propo:codeInclusion} is immediate.

\begin{coro}\label{coro:switchingCodes} Let $H'$ be obtained from a Hadamard matrix $H$ by switching a closed quadruple.
If the linear binary codes associated to the columns of $H$ and $H'$ are both self-dual, then they are equal.  In particular, if $H$ is of size $n\equiv8\pmod{16}$, then  the linear binary codes associated to the columns of $H$ and $H'$ are equal.
\end{coro}

Corollary~\ref{coro:switchingCodes} will be important when we discuss classification of
Hadamard matrices of size 24.  Since in the setting  of Corollary~\ref{coro:switchingCodes} there is a one-to-one correspondence between closed quadruples and weight 4 code
words, the corollary also provides an alternative
proof of Proposition~\ref{propo:CQinvariant}.

\subsection{A Hall set switching invariant}
\label{subsection:smith}
An important notion used in the classification of Hadamard matrices is that of
integer equivalence.

\begin{defin}
Two integer matrices $A$ and $B$ are {\em integer equivalent} if $A$ can be converted
to $B$ by some sequence of the following row and column operations:
\begin{itemize}
\item permutation of rows (columns)
\item negation of rows (columns)
\item addition of an integer multiple of a row (column) to another row (column).
\end{itemize}
\end{defin}

Associated to the integer equivalence class of a matrix $A$ of size $n$ is a set of integers $s_1,\ldots,s_n$
called {\em invariant factors} satisfying:
\begin{enumerate}
\item The matrix $\diag(s_1,\ldots,s_n)$ is integer equivalent to $A$.
\item There exists $r$ such that $1\le r\le n$ and $s_i|s_{i+1}$ for $1\le i\le r-1$ and $s_{r+1}=\ldots=s_n=0$.
\item The product $s_1s_2\ldots s_i$ equals the GCD of the $i\times i$ minors of $A$.
\end{enumerate}
The matrix $\diag(s_1,\ldots,s_n)$ is called the {\em Smith normal form} of $A$.  Two
integer equivalent matrices have the same Smith normal form.

A number of properties of the Smith normal form of a Hadamard matrix  have been
proved~\cite{Wal71,New71}:
\begin{enumerate}
\item $s_1=1$; $s_2=\ldots=s_{\alpha+1}=2$, for some $\alpha\ge\lfloor\log_2 n\rfloor+1$;
\item $s_i s_{n+1-i}=n$.
\end{enumerate}

In order 36, for example, we have~\cite{CooMilWal78}
\begin{itemize}
\item $s_1=1$
\item $s_i=2$ for the next $\alpha$ values of $i$.  ($2\le i \le \alpha+1$)
\item $s_i=6$ for the next $34-2\alpha$ values of $i$
\item $s_i=18$ for the next $\alpha$ values of $i$
\item $s_{36}=36$
\end{itemize}
where $6\le\alpha\le17$.  The single parameter $\alpha$ determines the integer equivalence
class of a Hadamard matrix $H$ in order 36, and we say that $H$ is in {\em Smith class} $\alpha$.

That the Smith class is invariant under switching Hall sets is implied by the following:

\begin{propo}\label{propo:smith}
If $B$ is obtained from $A$ by switching a Hall set, then $B$ is integer equivalent to $A$.
\end{propo}

\begin{proof}
Switching a Hall set can be achieved by a sequence of integer
row and column operations.  Let the order of the matrix in~(\ref{hallSet}) be $4k+4$.  Adding
each of rows 1 through 4 to each of the $k$ rows 5 through $k+4$, and then adding each
of columns 1 through 4 to each of columns 5 through $k+4$ has the effect of negating $F_1$ and
$G_1$.
\end{proof}

\section{Equivalence relations}
Hadamard equivalence, usually simply called ``equivalence,'' was defined in the introduction.  
We will call Hadamard equivalence classes {\em H-classes.}  By adjoining
additional operations to the list of operations given there, we can define new equivalence
relations.  We already did this in the previous section when we defined integer equivalence, whose
equivalence classes are the Smith classes.  The considerations of the previous section also allow us
to define {\em equivalence with respect to the associated binary linear code,} or {\em code equivalence} for short: Two Hadamard matrices are code equivalent if their binary linear codes, as defined in the statement of Proposition~\ref{propo:wordsToQuads} are isomorphic.  In this section,  we define further notions of
equivalence.

\begin{defin}
If $n\equiv0\pmod{8}$ then two Hadamard matrices $A$ and $B$ of size $n$ are {\em Q-equivalent}
if $B$ can be obtained from $A$ by some sequence of the operations
\begin{itemize}
\item row or column negation
\item row or column permutation
\item switching a closed quadruple of rows
\item switching a closed quadruple of columns.
\end{itemize}
If the last operation is disallowed, then $A$ and $B$ are said to be QR-equivalent; if the third operation
is disallowed then $A$ and $B$ are said to be QC-equivalent.
When $n\equiv4\pmod{8}$, Q-equivalence is defined by replacing the last two operations
with
\begin{itemize}
\item switching a Hall set.
\end{itemize}
Associated with these equivalence relations are equivalence classes, called {\em Q-classes,}
{\em QR-classes,} and {\em QC-classes.}
\end{defin}

Hadamard equivalence is stronger than Q-equivalence and therefore has a more
refined equivalence class structure.  In other words, there are at least as many H-classes as
there are Q-classes, and each H-class is contained entirely within a particular Q-class.
QR-equivalence (or QC-equivalence) is intermediate in strength
between H-equivalence and Q-equivalence, and will therefore have an intermediate
number of equivalence classes.  When $n\equiv8\pmod{16}$, Corollary~\ref{coro:switchingCodes}
implies that QR-equivalence is a refinement of code equivalence: two Hadamard matrices in the same
QR class have isomorphic codes; the converse does not necessarily hold as will be seen
in the case $n=24$, which is discussed in the next section.

By Proposition~\ref{propo:smith}, Q-equivalence is stronger than integer equivalence
when $n\equiv4\pmod{8}$ which implies that there are at least as many Q-classes as
there are Smith classes in those orders.

%Doubly-even binary codes can  be associated with Hadamard matrices.
%The following observations are due to Jennifer Key:  closed row quadruples correspond
%to  code words of weight 4 in the associated code;  in the case $n=24$
%the code is uniquely determined (up to equivalence) by the number of code words
%of weight 4.   When $n\equiv8\pmod{16}$ Proposition~\ref{propo:CQinvariant} implies that
%this number is an invariant under switching a closed row quadruple.  Hence the associated code is not 
%changed by switching.  What happens for $n=24$ is most
%likely a general phenomenon.  If so, then when $n\equiv8\pmod{16}$ the QR-classes will be
%a refinement of the equivalence classes associated to the doubly even binary code.

An equivalence class, of any type, may or may not be {\em self-dual.}
The {\em dual} of a set of matrices is the set containing their transposes.  A set that equals
its own dual is self-dual.
Many but not all Q-classes turn out to be self-dual.  In other words, many matrices
are Q-equivalent to their transposes.  From the row-column
symmetry in the definition of Q-equivalence it follows that if a Q-class contains at least one
self-dual matrix, then that Q-class is self-dual. 

We will see examples of these phenomena in the next section.

\section{Application to the enumeration of inequivalent Hadamard matrices}
We remind the reader that Hadamard matrices have been completely classified up to order
28.  There are five H-classes in order 16~\cite{Hal61}, three in order 20~\cite{Hal65}, 60 in order 24~\cite{ItoLeoLon81,Kim89}, and 487 in order 28~\cite{Kim94b}.  Using the available
lists of H-classes, which can be obtained from a number of sources~\cite{Seb05,Slo05,Spe05}, we will be able to determine the structure of the Q-classes in these orders.
The classification of H-classes in orders 32 and higher appears to be very difficult.  We will content
ourselves with identifying the Q-classes of all Hadamard matrices in orders 32 and 36 that were known
before the recent work of Bouyukliev, Fack, and Winne (see Introduction),
and completely enumerating those Q-classes that are small enough for this to be feasible.

Our procedure requires that we maintain a database of inequivalent matrices.  As new matrices
are generated, they are put in a canonical form
and compared with known matrices to prevent duplication in the database.
To put the matrices in canonical form, we followed the suggestion of Brendan
McKay~\cite{McK79}, converting $n\times n$ matrices to graphs on $4n$ vertices and then
using the graph isomorphism program {\em nauty} that he developed~\cite{McK81}.  The canonical
form of the graph computed by {\em nauty} was then converted back into a matrix.  As
suggested in the {\em nauty} User's Guide~\cite{McK04}, we used the vertex invariant {\em cellquads}
at level 2, which improves the efficiency in processing this type of graph.

To generate lists of inequivalent Hadamard matrices of order $n$ we carried out the following procedure,
which requires a seed Hadamard matrix of order $n$ as input:
\begin{enumerate}
\item Initialize {\tt hadList} to null list.
\item Compute canonical form of seed matrix using {\em nauty}.  Append it to {\tt hadList}.
\item Compute canonical form of transpose of seed matrix.  If it differs from canonical form
of seed matrix, append it to {\tt hadList}.
\item Initialize {\tt ctr} to 1.
\item Let $H$ be matrix number {\tt ctr} on {\tt hadList}.  If $n\equiv4\pmod{8}$ and this matrix is in the
H-class of the transpose of the previous one, skip to Step 7.
\item For each closed row quadruple ($n\equiv0\pmod{8}$) or Hall set ($n\equiv4\pmod{8}$) in $H$, 
\begin{enumerate}
\item Switch the quadruple (Hall set) and compute the canonical form of the resulting matrix to obtain $H'$.
\item If $H'$ differs from all matrices on {\tt hadList}, append it to {\tt hadList}.  Then if
the canonical form of the transpose of $H'$ differs from $H'$, append it to {\tt hadList} as well.
\end{enumerate}
\item Increment {\tt ctr}.  If {\tt hadList} is not exhausted, return to Step 5.
\end{enumerate}

Note that this procedure generates the Q-class of the seed matrix unless the Q-class
happens to be non-self-dual, in which case it generates the union of the Q-class and its dual.
This is due to the use of the transposition operation in Step 6(b).  
Non-self-dual Q-classes always turn out to be small, and when the situation arises,
we partition the union into two Q-classes by hand.  (We could use column
quadruple switching in the $n\equiv0\pmod{8}$ case and dispense with transposition in both cases,
thereby avoiding this issue, but we found it
convenient to use transposition to keep track of duality.)  We can also modify the procedure by simply eliminating the transposition step, in which case
the procedure generates the QR-class of the seed matrix in the $n\equiv0\pmod{8}$ case.

Here are the results on the Q-classes and QR-classes for orders $16$ and $24$:
\begin{itemize}
\item $n=16$: The five H-classes are all Q-equivalent.  More strikingly, they are all
QR-equivalent.
\item $n=24$: Of the 60 H-classes, 59 are Q-equivalent.  The H-class missing from the main
Q-class is that of the Paley matrix which has no closed quadruples and is self-dual.  It forms
a Q-class all by itself.

As stated in Section~\ref{subsection:codes}, Assmus and Key classified the 60 H-classes according to the doubly-even binary codes associated to the columns of the matrices.
(See Table 1 in~\cite{AssKey92a} or Table 7.1 in~\cite{AssKey92b}, but beware that
$42^{32}_D$, listed with the code $D$, should be listed with the code $C$, and that $32^{42}_D$ in
line 3 of the table should be changed to $32^{42}_C$.)  We use QR-equivalence to
refine this classification.

Assmus and Key found that six codes, labeled $A$, $C$, $D$, $E$, $F$, and $G$, occur.  They are distinguished
by the number of code words of weight 4 and therefore, according to Propositions~\ref{propo:wordsToQuads}
and~\ref{propo:quadsToWords}, by the number of closed row quadruples
in the associated matrices.  We now see Corollary~\ref{coro:switchingCodes}, on the invariance
of codes under switching, in action.  For example,
the matrices associated with the code $D$ all have 12 closed row quadruples.  Switching any
of these quadruples produces another matrix with code $D$.  Depending on which of these
matrices one starts with, switching row quadruples produces a QR-class of size 5 or of size 10.  These two
QR-classes together account for all 15 H-classes associated with the code $D$.  Results
for all the codes appear in Table~\ref{table:codes}.

Note that the matrices associated with the 
$[24,12]$ extended Golay code $G$ do not contain closed row quadruples.  One class of such matrices
must be that of the Paley Hadamard matrix as we have already stated that it has no closed
quadruples.  There is a second class of matrices with no closed row quadruples.  The matrices in
this class, however, do each have 66 closed column quadruples, since their duals turn out to
be in the class of the code $E$.
\end{itemize}

\begin{table}
\begin{tabular}{c|c|c|c}
 & \# weight-4 & size of & sizes of \\
code & code words & code class & QR-classes \\
\hline
A & 30 & 8 & 8 \\
C & 18 & 17 & 17 \\
D & 12 & 15 & 5, 10 \\
E & 66 & 8 & 8 \\
F & 6 & 10 & 5, 5 \\
G & 0 & 2 & 1, 1 \\
\end{tabular}

\vspace{6pt}

\caption{The 6 codes associated with the 60 Hadamard matrices of order 24.}
\label{table:codes}
\end{table}

The results on Q-classes in orders $20$ and $28$ are:
\begin{itemize}
\item $n=20$: The 3 H-classes are Q-equivalent.
\item $n=28$: Of the 487 H-classes, 486 of them (the ones containing Hall sets~\cite{Kim94a})
are Q-equivalent.   The Paley matrix (generated from quadratic residues in $\GF(3^3)$) contains
no Hall set and therefore its H-class forms a Q-class by itself.
\end{itemize}

Before presenting our results in orders 32 and 36, we ask what might the results so far
lead us to expect in higher order?  It is striking that except for a small number of exceptions
(the H-classes of the Paley
matrices in orders 24 and 28), all Hadamard matrices of given order are Q-equivalent.
Could this be a general phenomenon?

In order 36, a difficulty arises.  By Proposition~\ref{propo:smith} the Smith class is invariant
under the defining operations of Q-equivalence.  We will see that at least six different Smith
classes occur, and so there must be at least six Q-classes, each possibly containing
many H-classes.  The reason the multiplicity of Smith classes was not an issue in
order 28 is that $7=28/4$ is an odd square free number.
By a result in~\cite{WalWal69} this implies that all Hadamard matrices in order 28 lie in a single
Smith class.
From the foregoing discussion, the best we can hope for for general $n\equiv4\pmod{8}$
is that within each Smith class there
will be a single dominant Q-class, and that the total number of
Q-classes will still be small.

The results we have obtained so far  appear to support
the idea of a single dominant Q-class in order 32, and of a single dominant Q-class within each
Smith class in order 36.  The total number of Q-classes also appears to be very small relative to the
number of H-classes.  We found only a few tens of Q-classes in our analysis of the matrices known
prior to the work of Bouyukliev, Fack, and Winne, but a preliminary analysis of their matrices suggests
that the number will rise into the hundreds, if not higher.  Our method was
to collect as many Hadamard matrices as possible from the literature or using known
construction techniques, and then to apply our algorithm to each of these matrices in
order to obtain its Q-class.  In fortunate cases our program terminated in a reasonable
time, giving us a complete enumeration of the elements of the Q-class of the given seed
matrix. In less fortunate cases---and if our speculations are correct, this is expected to be the usual
situation---the Q-class was too big to enumerate completely.  Instead, we compared partially
constructed Q-classes with each other, and looked for overlaps.  By so doing, we managed
to identify unambiguously the Q-class of every Hadamard matrix in orders 32 and 36 known to us
prior to the work of Bouyukliev, Fack, and Winne, to enumerate
the smaller of these Q-classes, and to obtain lower bounds on the sizes of the larger Q-classes.

\subsection{Order 32}
\begin{propo}\label{propo:tensor32}
All Hadamard matrices  of either of the forms
\begin{equation*}
H=\begin{bmatrix}A & B\\A & -B\end{bmatrix}, \qquad
{\widetilde H}=\begin{bmatrix}A & A\\B & -B\end{bmatrix},
\end{equation*}
where $A$ and $B$ are any Hadamard matrices of order 16, are Q-equivalent.
\end{propo}
\begin{proof}
From the discussion in Section~\ref{subsection:tensor} it follows that, from the matrix
$\begin{bmatrix}A & B\\A & -B\end{bmatrix}$,
with $A$ and $B$ fixed, we may obtain any matrix of the form
$\begin{bmatrix}A & PB\\A & -PB\end{bmatrix}$,
by switching closed row quadruples.

To show that all matrices of the form $H$ are Q-equivalent we need only to show that
we can change the H-class of $A$ or of $B$ to any of the five classes in order 16 by switching closed quadruples.  Since
all Hadamard matrices of order 16 are QR- and QC-equivalent, we can achieve this by switching
closed column quadruples in the $A$ columns of $H$ only or in the $B$ columns of $H$ only.
(Closed column quadruples of $A$ or of $B$ extend to closed column quadruples of $H$
and switching transforms the top and bottom halves of a column the same way.)

Analogous arguments, with rows and columns interchanged, show that all matrices of the form
$\widetilde H$ are Q-equivalent.  To show that matrices of the form $H$ and of the form $\widetilde H$
are Q-equivalent to each other, simply note that both sets contain the Sylvester Hadamard matrix.
\end{proof}

Thus the 66099 H-classes identified in~\cite{LinWalLie93} are Q-equivalent.  We call
the Q-class of these matrices the {\em Sylvester Q-class}.  We now turn to other known
Hadamard matrices in order 32:
\begin{itemize}
\item the Paley matrix,
\item 13 matrices from generalized Legendre (GL) pairs~\cite{FleGysSeb01},
\item four matrices listed in~\cite{AraHarKha04} and their transposes,
\item the maximal excess matrix in~\cite{EnoMiy80},
\item four matrices from Construction II in~\cite{LinWalLie93},
\item two Williamson matrices,
\item eight Goethals-Seidel matrices constructed from circulant blocks,
\item 18 Goethals-Seidel matrices constructed from negacyclic blocks,
\item 10 matrices constructed from two circulants,
\item 17 matrices constructed from two negacyclic matrices,
\item a matrix from the appendix of~\cite{LamLamTon01} and its transpose.
\end{itemize}
Some of these matrices were provided by Hadi Kharaghani.
Discarding duplicates (which occur due to accidental equivalences) and matrices that happen to
have one of the forms
in Proposition~\ref{propo:tensor32}, we are left with a list of 59 matrices.  Of these, 49 are in
the Sylvester Q-class.  Using these matrices, and some matrices from Proposition~\ref{propo:tensor32}
as seeds, we have managed to generate 3,577,996 H-classes in the Sylvester Q-class by using our
program and then piecing together the results.  This
is certainly a gross underestimate of the actual number.

The ten exceptional matrices among the 59 all lack closed
quadruples either in rows or in columns, and therefore form Q-classes by themselves.
Of the ten exceptional matrices,
six are constructed from GL pairs, and four are constructed from two negacyclic blocks.  The matrices
from GL pairs are listed on the web page~\cite{Seb05} as P12--P19 (with transposes of non-self-dual
matrices omitted).
The exceptional GL pair matrices are P13, P15 and its transpose, P17, and P19 and its transpose.
Matrix P17 is Hadamard equivalent to the Paley matrix.
Of the matrices constructed from two negacyclic blocks, the exceptional ones come in
two dual pairs.

The Sylvester Q-class and the ten singleton Q-classes give total of 11 known Q-classes in
order 32, containing at least 3,578,006 Hadamard equivalence classes.

\subsection{Order 36}
As noted above, in order 36 we must consider each Smith class separately.  Although Smith classes
$\alpha=6,7,\ldots,17$ are allowed, the only Smith classes known to be nonempty are
$\alpha=11,12,13,14,15,16, 17$.  

A complete summary of the seed matrices we compiled in order 36 follows:
\begin{itemize}
\item Ted Spence's 180 matrices related to regular 2-graphs (S1--S180)~\cite{Spe95,McKSpe01,Spe05},
\item the 24 matrices of Goethals-Seidel type classified by Spence and
Turyn (GS1--GS24)~\cite{Spe05},
\item the 11 matrices with automorphism of size 17 classified by Tonchev
(T1--T11)~\cite{Ton86},
\item the Bush-type Hadamard matrix found by Janko (B1)~\cite{Jan01},
\item a regular Hadamard matrix found by Jennifer Seberry and listed on
her web page (R1)~\cite{Seb05}, (She actually lists four, but two are duplicates, and two
are of Goethals-Seidel type.)
\item four Williamson Hadamard matrices (W1--W4), (There is a fifth, but it is
equivalent to one of Tonchev's.)
\item the $(35,17,8)$-difference set construction (D1),
\item seven matrices of the type defined by Whiteman (a Goethals-Seidel
array bordered by a Hall set) (Wh1--Wh7)~\cite{Whi76},
\item two block negacyclic Bush-type Hadamard matrices, the first given
in the paper of Janko and Kharaghani (NB1, NB2)~\cite{JanKha02},
\item a matrix in Smith class 11, found in the course of a (fruitless)
search for block circulant Bush-type matrices (O1),
\item a skew Bush-type Hadamard matrix found by Leif J{\o}rgensen and its transpose (J1, J2)~\cite{Jor05},
\item a matrix listed in the appendix of~\cite{LamLamTon01} (LLT1),
\item the first matrix known (to us) in Smith class 16, found by Bouyukliev, Fack, and Winne~\cite{BouFacWin05} (BFW1).
\end{itemize}
Reference~\cite{CooMilWal78}
was helpful in assembling the above list, but the reader should note that the 80 matrices from Steiner triple systems, which
are a subset of Spence's 180 matrices, are in Smith class 13, not 12 as stated there.  We have not made
a serious effort to credit the original author of every matrix on our list, as we were more
concerned with compiling as complete a list as possible from readily obtainable sources.  We should
note, however, that many of these matrices derive from the important work of Goethals and Seidel~\cite{GoeSei70}, including the 80 matrices from Steiner triple systems mentioned above, and 11 matrices derived from Latin squares of order 6, which are also a subset of Spence's list.

The structure we have uncovered in Smith class 13 is interesting, so we describe it in detail.  Hadamard matrices in this class include
179 of Ted Spence's 180 matrices.  (His matrix 137 is in Smith
class 11.)  Two other matrices in Smith class 13 were previously known: the regular Hadamard matrix
constructed by Seberry, and the block negacyclic Bush-type Hadamard matrix
constucted by Janko and Kharaghani.  Seberry's matrix and 172 of Spence's fall into
the same Q-class which we found has size 3425.  Two of Spence's matrices (179 and 180) and the Bush-type matrix
form singleton Q-classes.  They have no Hall sets.  The remaining five of Spence's matrices
lie in a Q-class of size 6.

Spence's matrix 137, which is in Smith class 11 and is one of the matrices derived from a Latin square of order 6, is intriguing.  It has nine Hall sets, but switching any
of these produces a matrix H-equivalent to the original.

Only two matrices on our list are in Smith class 14, B1 and LLT1.  They are Q-equivalent.
A major success of our program has been the complete enumeration of their Q-class, which has 954,254 elements.  In each of Smith classes 15, 16, and 17 there is one known Q-class of size above five million, while all
other known Q-classes are of size no greater than 5.  The three large Q-classes have not yet been
completely enumerated.  At present, there is no evidence for more than one large Q-class
in any Smith class.

The 236 matrices we compiled represent seven different Smith classes, and lie in 21
different Q-classes.  Some details are given in Table~\ref{table:order36}.  The union of
the known Q-classes contains at least 18,292,717 Hadamard equivalence classes of order 36.

%Note: After this work was completed, Bouyukliev, Fack, and Winne announced the complete classification
%of $2$-$(35,17,8)$ designs with an automorphism of order 3 fixing two points and blocks.  From these
%designs, they find 7238 H-equivalence classes~\cite{BouFacWin05}.  Most of these matrices have not yet been analyzed.  Such
%an analysis should provide a good test of the ideas of this  paper regarding using Q-equivalence
%in classifying Hadamard matrices.  Among the matrices of Bouyukliev, Fack, and Winne are the first  examples known (to this author) in Smith class 16.  One of these was used to generate
%a Q-class of size at least 1,010,890.  This brings the number of known Smith classes to 7, the
%number of known Q-classes to 21, and the number of known H-equivalence classes to 4,745,357.
%Preliminary analysis of a sampling of the the new matrices of Bouyukliev, Fack, and Winne indicates the presence of a number of
%new small Q-classes.  We intend to make a complete enumeration of these, and a full analysis of all the new matrices.  The results will be presented in a follow-up to the present paper.

\begin{table}
\begin{tabular}{l|l}
$\alpha$ & Q-classes \\
\hline
11 & 1 (S137), 1 (O1) \\
12 & 1 (D1) \\
13 & 1 (S179), 1 (S180), 1 (NB1), 6 (S172), 3425 (S1) \\
14 & $954,254$ (B1) \\
15 & 5 (W3), 5 (W4), $\ge5,520,880$ (GS1) \\
16 & $\ge5,814,129$ (BFW1) \\
17 & 1 (GS11), 1 (GS12), 1 (T1), 1 (T2), 1 (T5), 1 (T6), 1 (T7),
 $\ge6,000,000$ (GS4) \\
\end{tabular}

\vspace{6pt}

\caption{Sizes of known Q-classes in order 36 for the 6 known Smith
classes, $\alpha$.  A representative matrix is listed for each Q-class.}
\label{table:order36}
\end{table}

\section*{Acknowledgments}
I thank Bruce Solomon for numerous interesting discussions and for
careful comments on the manuscript.  I thank Hadi Kharaghani for extensive
correspondence and for providing many unpublished Hadamard matrices in orders 32 and 36.
I thank Jennifer Key for recomputing a table that appears in~\cite{AssKey92a,AssKey92b}, and for pointing out the connection between closed quadruples and code
words of weight 4.  I am indebted to Robin Chapman who pointed out an error in the
proof of Proposition~\ref{propo:hallSwitch} in an earlier version of the manuscript and who suggested the simple proof of Proposition~\ref{propo:closed}.  Iliya Bouyukliev, Veerle Fack, and Joost Winne kindly provided me with a list of the matrices they obtained in their work on classification of $2$-$(35,17,8)$-designs.  Jayoung Nam made numerous helpful suggestions about presentation.
I thank the editor, Kevin Phelps, and the anonymous referees for suggestions that improved the paper.
This project relied on the High Performance Computing facilities of
Indiana University, in particular the IBM RS/6000 SP system.

\bibliographystyle{plain}
\bibliography{hadgenVer4}

\end{document}